\documentclass[3p,times]{elsarticle}
\usepackage[figuresright]{rotating}
\usepackage{enumerate}

\usepackage{epstopdf}
\usepackage{caption}
\usepackage{subfigure}
\usepackage{amsmath}
\usepackage{dcolumn}
\usepackage{float}

\usepackage{amssymb}

\usepackage{subfigure}

\renewcommand{\figurename}{Fig.}

\usepackage{amsmath}

\usepackage{array}
\usepackage{amsfonts}
\usepackage{bm}
\usepackage{algorithm}
\usepackage{algorithmic}
\usepackage{amscd}
\usepackage{mathrsfs}

\usepackage{amssymb}
 \usepackage{amsthm}
 \usepackage{amsmath}
\usepackage{enumerate}
\usepackage[figuresright]{rotating}
 \usepackage[titletoc]{appendix}
 \usepackage{appendix}

\renewcommand{\figurename}{Fig.}

\newcommand{\PreserveBackslash}[1]{\let\temp=\\#1\let\\=\temp}
\newcolumntype{C}[1]{>{\PreserveBackslash\centering}p{#1}}
\newcolumntype{R}[1]{>{\PreserveBackslash\raggedleft}p{#1}}
\newcolumntype{L}[1]{>{\PreserveBackslash\raggedright}p{#1}}
\begin{document}

\begin{frontmatter}



%% e.g. \dochead{17th International Conference on Dynamical Processes in Excited States of Solids}

\title{Numerical methods for the two-dimensional Fokker-Planck equation governing the probability density function of the tempered fractional Brownian motion
%Variational properties and effective calculation for problems involving  the tempered fractional  Calculus %\footnote{This work was supported by the National Natural Science Foundation of China under Grant  No. 11271173 and No. 11471150.}
}


\author{Xing Liu and Weihua Deng}

\address{School of Mathematics and Statistics, Gansu Key Laboratory of Applied Mathematics and Complex Systems, Lanzhou University, Lanzhou 730000, P.R. China}

\begin{abstract}
In this paper, we study the numerical schemes for the two-dimensional Fokker-Planck equation governing the probability density function of the tempered fractional Brownian motion. The main challenges of the numerical schemes come from the singularity in the time direction. When $0<H<0.5$, a change of variables $\partial \left(t^{2H}\right)=2Ht^{2H-1}\partial t$ avoids the singularity of numerical computation at $t=0$, which naturally results in nonuniform time discretization and greatly improves the computational efficiency.
For $0.5<H<1$, the time span dependent numerical scheme and nonuniform time discretization are introduced to ensure the effectiveness of the calculation and the computational efficiency. By numerically solving the corresponding Fokker-Planck equation, we obtain the mean squared displacement of stochastic processes, which conforms to the characteristics of the tempered fractional Brownian motion.

%The tempered fractional operators are the generalization of the standard fractional Riemann-Liouville and Caputo operators, which include the smooth truncation effects to suit  the real physical problems.
%   In this paper, we first construct the    variational properties of the tempered  fractional derivatives and  integrals, and then  apply them to solve the  space and (or) time tempered fractional PDEs, including   the related theoretical analysis and effective numerical implementation.

\end{abstract}


\end{frontmatter}
%\makeatletter
%\def\@captype{figure}
%\makeatother
%

% ----------------------------------------------------------------
\section{Introduction}
Revealing the law of motion of particles in the system is always a hot subject due to its wide applications in physics, biology, chemistry, etc. The modeling of particles' motion can be traced back to Brownian motion, which is a normal diffusion process. With the advance of scientific research, more and more scientists realized that the motion of particles in complex disordered systems generally exhibits anomalous dynamics. The mean squared displacement (MSD) is usually used to distinguish the types of stochastic processes. The MSD of Brownian motion goes like $\left\langle (\Delta x)^2\right\rangle=\left\langle [x(t)-\langle
x(t)\rangle]^2\right\rangle\sim t^\nu$ with $\nu=1$; for a long time $t$ if $\nu\not=1$, it is called anomalous diffusion, being subdiffusion for $\nu<1$ and superdiffusion for $\nu>1$; in particular, it is termed as localization diffusion if $\nu=0$ and ballistic diffusion if $\nu=2$ \cite{1,3,4}. There are two types of typical stochastic processes to model anomalous diffusion: Gaussian processes and non-Gaussian ones \cite{2,3-4,5,6,7}. In order to obtain MSD, one must know the probability density function (PDF) of stochastic processes, which can be obtained by solving the corresponding Fokker-Planck equations.
%. Thus, numerical methods of the Fokker-Planck equation for stochastic processes are very important tool in the study of particles motion.
Studying the numerical methods of solving diffusion equations is a very active field, and many effective schemes have been developed \cite{8,9,10,11,12,13,14,15}. In this paper, we mainly focus on the numerical schemes for the newly developed models \cite{3}, i.e., the two-dimensional Fokker-Planck equations governing the PDF of the tempered fractional Brownian motion (tfBm), and reveal the mechanism of the motion of the particles.

The tfBm is defined as
\begin{displaymath}
B_{a,\lambda}(t)=\int_{-\infty}^{+\infty}\left[\mathrm{e}^{-\lambda(t-z)_+}(t-z)_+^{-\alpha}
-\mathrm{e}^{-\lambda(-z)_+}(-z)_+^{-\alpha}\right]B(\mathrm{d}z),
\end{displaymath}
describing anomalous diffusion with exponentially tempered long range correlations \cite{16,17}, where $B(z)$ is Brownian motion; and the tfBm is a Gaussian process with the corresponding Fokker-Planck equation \cite{3}
 \begin{displaymath}
 \frac{\partial P(x,t)}{\partial t}=\frac{\Gamma(H+1/2)}{\sqrt{\pi}(2\lambda)^H}\lambda t^HK_{H-1}(\lambda t)\frac{\partial^2P(x,t)}{\partial x^2},
 \end{displaymath}
 where
 \newcommand{\ud}{\mathrm{d}}
 \begin{equation}\label{eq:1}
 K_H(\lambda t)=\frac{1}{2}\int_0^{\infty} z^{H-1}\exp\left[-\frac{1}{2}\lambda t\left(z+\frac{1}{z}\right)\right]\ud z,
 \end{equation}
 the Hurst index $H=0.5-\alpha$, and $0<H<1, H\not=0.5, \lambda>0$. %The density distribution function $P(x,t)$ describes distribution of particles at time $t$.

 Along the direction of extension of the one-dimensional Fokker-Planck equation, we get
 \begin{equation}\label{eq:2}
 \frac{\partial u(x,y,t)}{\partial t}=\frac{\Gamma(H+1/2)}{\sqrt{\pi}(2\lambda)^H}\lambda t^HK_{H-1}(\lambda t)\left[\frac{\partial^2}{\partial x^2}+\frac{\partial^2}{\partial y^2}\right]u(x,y,t),
 \end{equation}
 with the initial and boundary conditions given by
 \begin{eqnarray*}
 u(x,y,0)=u_0(x,y),\quad(x,y)\in\Omega,\\
 u(x,y,t)=0,\quad (x,y,t)\in\partial\Omega\times[0,T].
 \end{eqnarray*}
 Here $\Omega=(0,L)\times(0,L')$ is the spatial domain, $\partial\Omega$ is the boundary of $\Omega$. The analytical solution of Eq.\ \eqref{eq:2} is difficult to find and one has to resort to numerical schemes. In the literatures \cite{18,19,20,21,22,23,24,25}, the explicit method, implicit method, and ADI method, etc have been discussed, however, the diffusion coefficients are usually constant.
 % Some researches on the numerical approximation of diffusion equation are known (explicit method, implicit method, ADI method, etc) \cite{18,19,20,21,22,23,24,25}. However, the diffusion coefficients of the diffusion equations are usually constant in these researches.
 In Eq.\ \eqref{eq:2}, the diffusion coefficient tends to infinity as $t\to 0$ and $0<H<0.5$; and the diffusion coefficient increases first and then decreases with time $t$ if $0.5<H<1$. For these two reasons, it is difficult to solve the equation directly using classical methods. We hope to obtain effective numerical methods by analyzing and applying the properties of diffusion coefficients.

 This paper is organized as follows. In section 2, for $0<H<0.5$, we derive a modified implicit method to circumvent the singularity of numerical calculation at $t=0$ and use nonuniform time stepsizes to improve computational efficiency and ensure the accuracy of the numerical solution. As $0.5<H<1$, a time span dependent numerical method with nonuniform time stepsizes is proposed to improve the accuracy of numerical solution and computational efficiency. In section 3, we present numerical results and show the characteristics of diffusion process corresponding to Eq.\ \eqref{eq:2}. Finally, a brief conclusion is provided in section 4.
\section{Numerical schemes with nonuniform time stepsizes}
Now, we introduce the numerical schemes that not only eliminate the singularity, but also improve the computational efficiency and ensure the accuracy of the numerical solution.
Generally in designing numerical schemes, one first needs to get a mesh in the space-time region where one wants to acquire the numerical approximation $u_{m,n}^k$ of the exact solution $u\left(x_m,y_n,t_k\right)$, where $\left(x_m,y_n,t_k\right)$ is the coordinate of the $(m,n,k)$ node of the mesh, and
$U_{M,N}^K=\left\{u_{m,n}^k,\, 0\le m\le M,\, 0\le n\le N,\, 0\le k\le K\right\}$.  In order to facilitate the numerical calculations, we rewrite the matrix $U_{M,N}^k$ as a vector $\widetilde{U}_{M,N}^k=\left(u_{1,1}^k\cdots u_{m,1}^k, u_{1,2}^k\cdots u_{m,2}^k\cdots  u_{1,n}^k\cdots u_{m,n}^k\right)^\mathrm{T}$.

Since in space the solution of Eq. \eqref{eq:2} has homogeneous properties,
%Because the direction of space does not cause trouble for numerical calculations,
the sizes of the mesh $\Delta x=x_{m+1}-x_m=h$ and $\Delta y=y_{n+1}-y_n=l$ are taken as  constants, and we discretize the operators $\frac{\partial^2}{\partial x^2}$ and $\frac{\partial^2}{\partial y^2}$ by means of the three-point centered formulas, i.e.,
\begin{displaymath}
\delta_x^2u_{m,n}^k=u_{m+1,n}^k-2u_{m,n}^k+u_{m-1,n}^k,
\quad\delta_y^2u_{m,n}^k=u_{m,n+1}^k-2u_{m,n}^k+u_{m,n-1}^k.
 \end{displaymath}
In the following, we sufficiently make use of the properties of the time dependent diffusion coefficients to do the time discretizations.
% Since diffusion coefficient is a function that takes time as an independent variable, in order to make the numerical solution effective, we have to consider the properties of the diffusion coefficient to choose time step.
%By magnifying and shrinking $t^HK_H(\lambda t)$, we have
For $t^HK_H(\lambda t)$, there exist the estimates
 \begin{eqnarray*}
t^HK_H(\lambda t)&=&\frac{t^H}{2}\int_0^\infty z^{H-1}\exp\left[-\frac{1}{2}\lambda t\left(z+\frac{1}{z}\right)\right]\ud z \\
&\le&\frac{t^H}{2}\int_0^\infty z^{H-1}\exp\left[-\frac{1}{2}\lambda tz\right]\ud z\\
&=&\frac{t^H}{2}\int_0^\infty \left(\frac{1}{2}\lambda tz\right)^{H-1}\left(\frac{1}{2}\lambda t\right)^{1-H}\exp\left[-\frac{1}{2}\lambda tz\right]\ud z\\
&=&\frac{2^{H-1}}{\lambda^H}\Gamma(H)
\end{eqnarray*}
and
\begin{eqnarray*}
t^HK_H(\lambda t)&\ge&\frac{t^H}{2}\int_1^\infty z^{H-1}\exp\left[-\frac{1}{2}\lambda t\left(z+\frac{1}{z}\right)\right]\ud z \\
&\ge&\frac{t^H}{2}\int_1^\infty z^{H-1}\exp\left[-\frac{1}{2}\lambda t(z+1)\right]\ud z \\
&=&\frac{2^{H-1}\mathrm{e}^{-\frac{1}{2}\lambda t}}{\lambda^H}\Gamma\left(H,\frac{\lambda t}{2}\right).
\end{eqnarray*}
Consequently, we have
\begin{displaymath}
\frac{2^{H-1}\mathrm{e}^{-\frac{1}{2}\lambda t}}{\lambda^H}\Gamma\left(H,\frac{\lambda t}{2}\right)\le t^HK_H(\lambda t)\le\frac{2^{H-1}}{\lambda^H}\Gamma(H),
\end{displaymath}
which leads to
\begin{equation}\label{eq:3}
\lim_{t \rightarrow 0}t^HK_{H}(\lambda t)=\frac{2^{H-1}}{\lambda^H}\Gamma(H).
\end{equation}
 Eq.\ \eqref{eq:1} shows that
 \begin{equation}\label{eq:4}
 K_{H-1}(\lambda t)= K_{1-H}(\lambda t),
 \end{equation}
combining with Eq.\ \eqref{eq:3} %and \eqref{eq:4}
 results in
 \begin{displaymath}
 \lim_{t \rightarrow 0}t^{1-2H}\left(t^HK_{H-1}(\lambda t)\right)=\frac{2^{-H}}{\lambda^{1-H}}\Gamma(1-H).
 \end{displaymath}
 Therefore, we need to, respectively, design the difference schemes of Eq.\ \eqref{eq:2} in two different cases, i.e., $2H-1<0$ and $2H-1>0$.

 $Case\ I$: As $0<H<0.5$, $\lim\limits_{t \rightarrow 0}t^HK_{H-1}(\lambda t)$ diverges. In order to eliminate the singularity, multiplying both sides of Eq.\ \eqref{eq:2} by $t^{1-2H}$, we get
 \begin{equation}\label{eq:5}
 \frac{\partial u(x,y,t)}{\partial \left(t^{2H}\right)}=\frac{\Gamma(H+1/2)}{2H\sqrt{\pi}(2\lambda)^H}\lambda t^{1-H}K_{H-1}(\lambda t)\left[\frac{\partial^2}{\partial x^2}+\frac{\partial^2}{\partial y^2}\right]u(x,y,t).
 \end{equation}
In this case of Eq. \eqref{eq:2}, with the increase of the time, the diffusion coefficient decays approximately as power law at a small time and exponentially at a relatively large time. To balance the decay of diffusion coefficient and make the variation of the solution approximately stationary, by taking $t^{2H}$ as a whole variable,   %the equation gradient approximately stationary,
we get a nonuniform discretization of $[0,T]$ with $t_k=(\tau k)^{1/{2H}}, \tau>0$, which greatly reduces the computation cost while keeping the accuracy.

%. Therefore, the time division of the above nonuniform method can improve the accuracy of the numerical solution and reduce the number of iterations of the calculations.

 From now on, the finite difference scheme can be obtained by discretizing Eq.\ \eqref{eq:5}:
\begin{equation}\label{eq:6}
\frac{u_{m,n}^{k+1}-u_{m,n}^k}{\bigtriangleup \left(t_k^{2H}\right)}=\frac{\Gamma(H+1/2)}{2H\sqrt{\pi}(2\lambda)^H}\lambda t_{k+1}^{1-H}K_{H-1}(\lambda t_{k+1})\left[\frac{\delta_x^2}{(\bigtriangleup x)^2}+\frac{\delta_y^2}{(\bigtriangleup y)^2}\right]u_{m,n}^{k+1},
\end{equation}
 which can be rearranged as
\begin{equation}\label{eq:7}
\left(1+\frac{2r}{(\bigtriangleup x)^2}+\frac{2r}{(\bigtriangleup y)^2}\right)u_{m,n}^{k+1}-\frac{r}{(\bigtriangleup x)^2}\left(u_{m+1,n}^{k+1}+u_{m-1,n}^{k+1}\right)-\frac{r}{(\bigtriangleup y)^2}\left(u_{m,n+1}^{k+1}+u_{m,n-1}^{k+1}\right)=u_{m,n}^k,
\end{equation}
where
\begin{displaymath}
r=\frac{\Gamma(H+1/2)}{2H\sqrt{\pi}(2\lambda)^H}\left[\lambda t_{k+1}^{1-H}K_{H-1}(\lambda t_{k+1})\tau\right].
\end{displaymath}
The coupling form of Eq.\ \eqref{eq:7} can be written as
\begin{displaymath}
C(t_{k+1})\widetilde{U}_{M,N}^{k+1}=\widetilde{U}_{M,N}^k,
\end{displaymath}
implying  that
\begin{displaymath}
\widetilde{U}_{M,N}^{k+1}=\prod_{j=0}^kC^{-1}(t_{j+1})\widetilde{U}_{M,N}^0,
\end{displaymath}
where $C(t_{k+1})$ is a growth matrix, being  symmetrical. The specific form of $C(t_{k+1})$ is omitted for the sake of brevity.

$Case\ II$: As $0.5<H<1$, the function $\lambda^{1-H}t^HK_{H-1}(\lambda t)$ increases first and then decreases with time. We define
\begin{displaymath}
t_{\max}:=\max\limits_{t_* \in [0,T] }\left\{\lambda^{1-H}t_*^HK_{H-1}(\lambda t_*) \ge\lambda^{1-H}t^HK_{H-1}(\lambda t) \quad {\rm for~any} \,\, t\in[0,t_*]\right\}.
\end{displaymath}
The maximum point of the function $\lambda^{1-H}t^HK_{H-1}(\lambda t)$ is (if it is not $T$)
%The function $\lambda^{1-H}t^HK_{H-1}(\lambda t)$ has an extreme point at time $t_{\max}$ which has been estimated using a large number of $H$ and $\lambda$, and then fit to the following equation.
\begin{equation}
t_{\max}\approx\frac{0.7442H-0.148H^{-1.3075}}{\lambda}.
\end{equation}
In the interval $t\in[0,t_{\max}]$, the diffusion coefficient increases approximately as power law,  while in the interval $t\in[t_{\max},T]$, the diffusion coefficient decays exponentially. Thus, in order to balance the trend of diffusion coefficient and improve the accuracy of the numerical solution and computational efficiency, we introduce the time span dependent difference schemes to solve equation. We choose a nonuniform partitions of $[0,t_{\max}]$ with $t_k=(\tau k)^{1/{2H}}$ being the power law decay and a nonuniform partition of $[t_{\max},T]$ with $t_k=(\tau k)^{1/{H}}$,  which is power law increasing.
%The numerical scheme of time span dependent nonuniform time step is more effective than other methods, and the number of the numerical iterations is cut down.

One can set up the difference scheme of Eq.\ \eqref{eq:2} as
%Rewriting Eq.\ \eqref{eq:2} to cooperate with two steps nonuniform timesteps method, and letting $t_{\max}=t_{k_1}$, one can set up the difference scheme
\begin{equation}\label{eq:9}
\left\{
\begin{array}{l}
\frac{u_{m,n}^{k+1}-u_{m,n}^k}{\bigtriangleup \left(t_k^{2H}\right)}=\frac{\Gamma(H+1/2)}{2H\sqrt{\pi}(2\lambda)^H}\lambda t_{k+1}^{1-H}K_{H-1}(\lambda t_{k+1})\left[\frac{\delta_x^2}{(\bigtriangleup x)^2}+\frac{\delta_y^2}{(\bigtriangleup y)^2}\right]u_{m,n}^{k+1},\quad k\le k_1,\\
\frac{u_{m,n}^{k+1}-u_{m,n}^k}{\bigtriangleup \left(t_k^{H}\right)}=\frac{\Gamma(H+1/2)}{H\sqrt{\pi}(2\lambda)^H}\lambda t_{k+1}K_{H-1}(\lambda t_{k+1})\left[\frac{\delta_x^2}{(\bigtriangleup x)^2}+\frac{\delta_y^2}{(\bigtriangleup y)^2}\right]u_{m,n}^{k+1},\quad k\ge k_1 .
\end{array}
\right.
\end{equation}
The coupled form of \eqref{eq:9} can be written as, when $k\le k_1$,
\begin{displaymath}
C_1(t_{k+1})\widetilde{U}_{M,N}^{k+1}=\widetilde{U}_{M,N}^k,\quad
\widetilde{U}_{M,N}^{k+1}=\prod_{j=0}^kC_1^{-1}(t_{j+1})\widetilde{U}_{M,N}^0,
\end{displaymath}
and when $k\ge k_1$,
\begin{displaymath}
C_2(t_{k+1})\widetilde{U}_{M,N}^{k+1}=\widetilde{U}_{M,N}^k,\quad
\widetilde{U}_{M,N}^{k+1}=\prod_{j_1=k_1}^kC_2^{-1}(t_{j_1+1})\prod_{j=0}^{k_1-1}C_1^{-1}(t_{j+1})\widetilde{U}_{M,N}^0.
\end{displaymath}
Growth matrices $C_1(t_{k+1})$ and $C_2(t_{k+1})$ correspond to two equations in \eqref{eq:9}.

By solving \eqref{eq:6} or \eqref{eq:9}, one can get all $u_{m,n}^k$, the approximations of the exact solution.
%Through the numerical calculations of \eqref{eq:6} and \eqref{eq:9}, we know all the $u_{m,n}^k$.
Checking stability and convergence is critical to understand the effectiveness of the numerical schemes. We use Fourier method to analyze the stability and convergence of the schemes \eqref{eq:6} and \eqref{eq:9} (for the details, see Appendix).
Let $e_{m,n}^k=u_{m,n}^k-u(x_m,y_n,t_k)$ be the difference between the numerical solution and the exact solution. The local truncation error is
\begin{equation}\label{eq:10}
R_{m,n}^{k}=O\left(\tau+h^2+l^2\right).
\end{equation}
We prove that the schemes \eqref{eq:6} and \eqref{eq:9} are unconditionally stable, that is,
\begin{equation}\label{eq:11}
\left\|U^k(x,y)\right\|_{L^2}^2<\left\|U^0(x,y)\right\|_{L^2}^2,
\end{equation}
and have first-order convergence in time and second-order convergence in space, i.e.,
\begin{equation}\label{eq:12}
\left\|e^k(x,y)\right\|_{L^2}^2<O\left(\tau+h^2+l^2\right).
\end{equation}

\captionsetup[figure]{name={Fig.}}
\begin{figure}[H]
\centering
\includegraphics[scale=0.5]{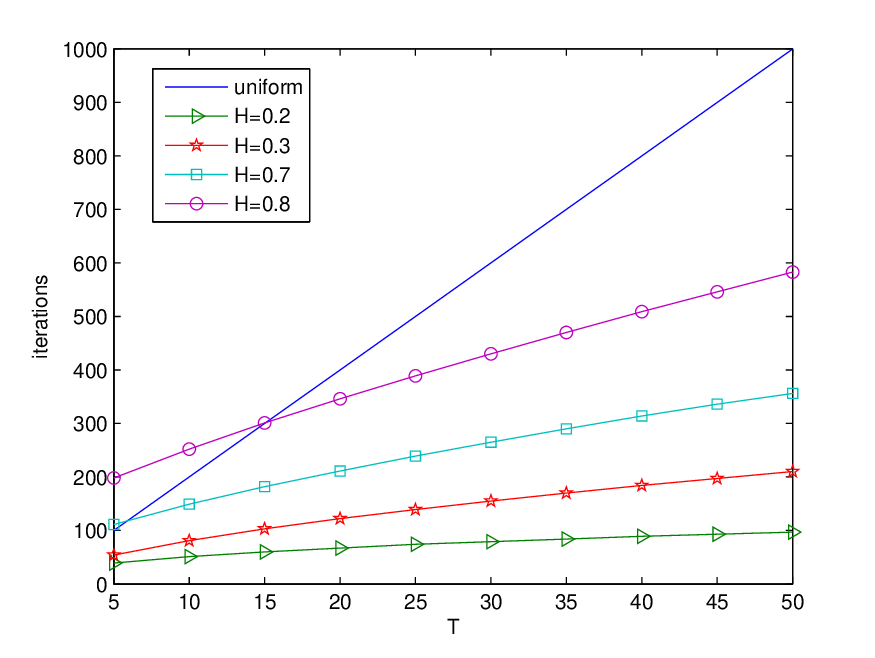}
\caption{Number of iterations (steps) for uniform and nonuniform ($H=0.2, 0.3, 0.7, 0.8$) time stepsizes with $\tau=0.05$ and $\lambda=0.1$.}\label{fig:1}
\end{figure}

\section{Localization diffusion: numerical results}
In this section, we test the schemes \eqref{eq:6} and \eqref{eq:9} by solving Eq. \eqref{eq:2} and calculating the MSD of tfBm. Although particles diffuse in unbounded domain, it is of course impossible to use boundary conditions at infinity in numerical calculations. Thus, we solve \eqref{eq:2} in a large enough two-dimensional domain $\Omega=(-100,100)\times(-100,100)$ by the proposed method. The numerical results $U_{M,N}^K$ with an initial data $u(x,y,0)=\mathrm{e}^{-x^2-2y^2}$ are obtained.

From \figurename\ \ref{fig:1}, we can see that the nonuniform time stepsize method is more computationally efficient. In order to study the motion of the particles, it is necessary to know the MSD of diffusion process tfBm. We define the MSD of two-dimensional stochastic process as
\begin{displaymath}
\left\langle \left[x(t)-\langle x(t)\rangle\right]^2+\left[y(t)-\langle y(t)\rangle\right]^2\right\rangle.
\end{displaymath}
The MSD of the stochastic process can be obtained by the numerical solution $U_{M,N}^K$ . After normalizing $U_{M,N}^k$, the discrete probability distribution $\mathrm{Pr}_{M,N}^K=\left\{\mathrm{Pr}_{m,n}^k, 0\le m\le M, 0\le n\le N, 0\le k\le K\right\}$ of particles is captured at each moment. The expectation formula implies that
\begin{displaymath}
\left\langle x(t_k)\right\rangle=\sum_{m=0}^M\sum_{n=0}^Nx_m\mathrm{Pr}_{m,n}^k,\quad
\left\langle y(t_k)\right\rangle=\sum_{n=0}^N\sum_{m=0}^My_n\mathrm{Pr}_{m,n}^k,
\end{displaymath}
and
\begin{displaymath}
\left\langle x^2(t_k)\right\rangle=\sum_{m=0}^M\sum_{n=0}^Nx_m^2\mathrm{Pr}_{m,n}^k,\quad
\left\langle y^2(t_k)\right\rangle=\sum_{n=0}^N\sum_{m=0}^My_n^2\mathrm{Pr}_{m,n}^k,
\end{displaymath}

\begin{figure}[H]
\centering
\subfigure[]{
\begin{minipage}[t]{0.5\textwidth}\label{fig:a}
\centering
\includegraphics[scale=0.4]{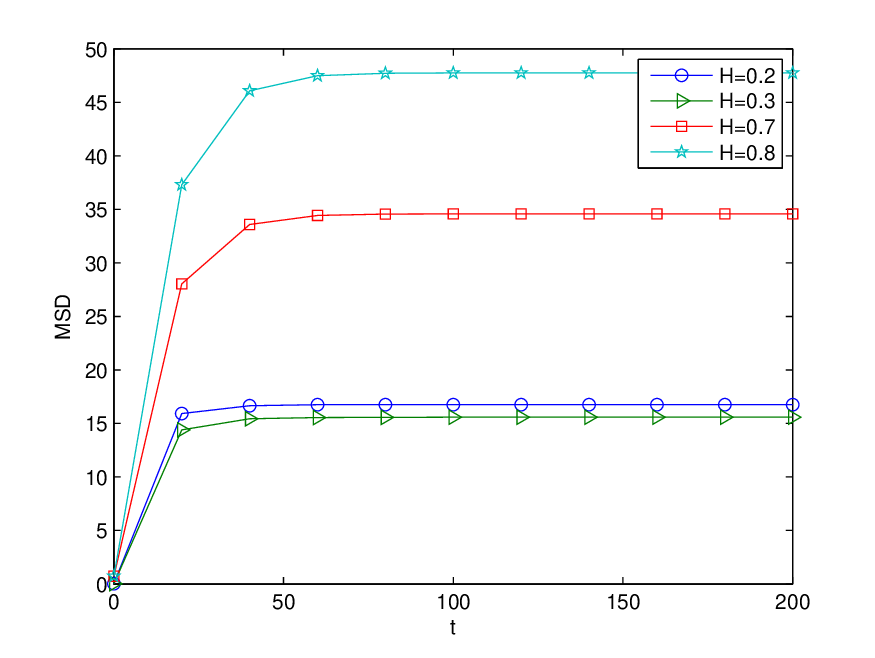}
\end{minipage}}\subfigure[]{
\begin{minipage}[t]{0.5\textwidth}\label{fig:b}
\centering
\includegraphics[scale=0.4]{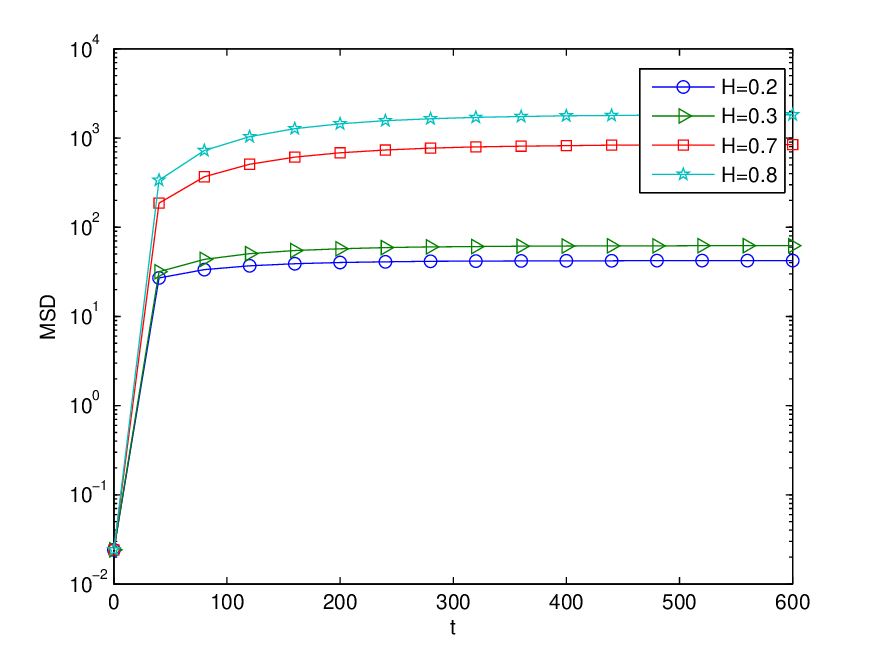}
\end{minipage}}
\caption{Simulations of MSD sampled over 1000 trajectories with $T=200$, $\lambda=0.1$ (left), $T=600$ and $\lambda=0.01$ (right).}\label{fig:2}
\end{figure}
\begin{figure}[H]
\centering
\subfigure[]{
\begin{minipage}[t]{0.5\textwidth}\label{fig:a2}
\centering
\includegraphics[scale=0.4]{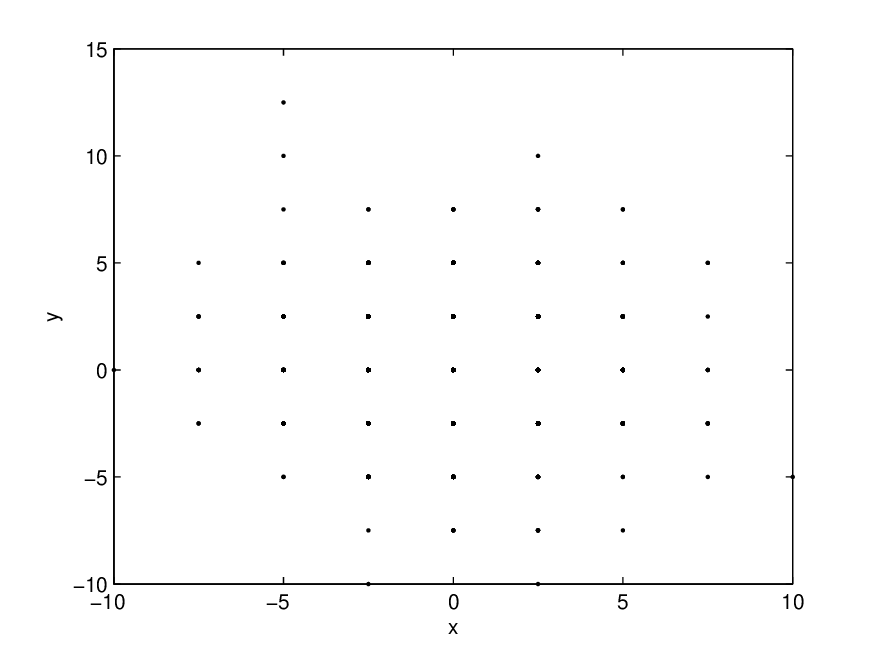}
\end{minipage}}\subfigure[]{
\begin{minipage}[t]{0.5\textwidth}\label{fig:b2}
\centering
\includegraphics[scale=0.4]{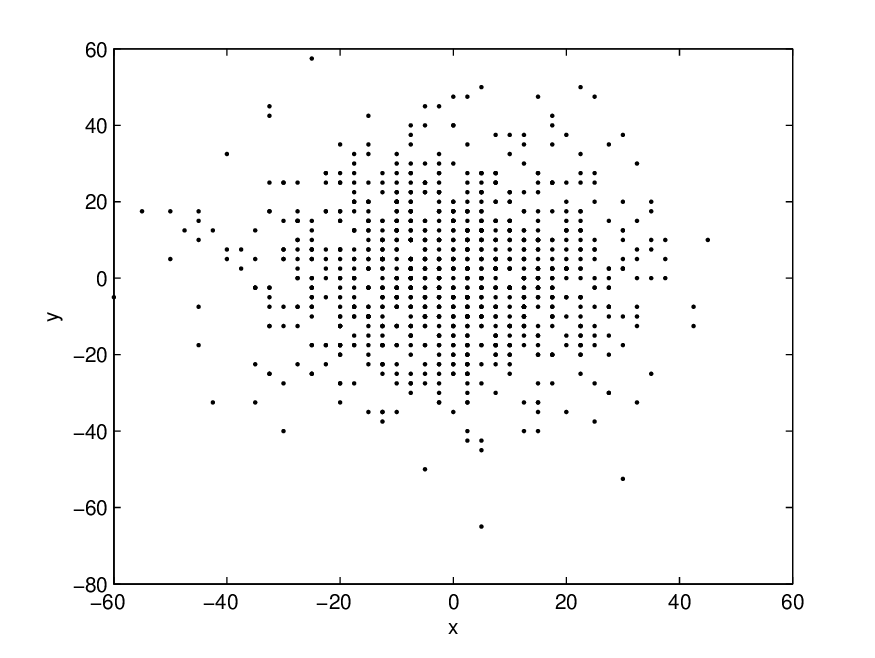}
\end{minipage}}
\subfigure[]{
\begin{minipage}[t]{0.5\textwidth}\label{fig:c2}
\centering
\includegraphics[scale=0.4]{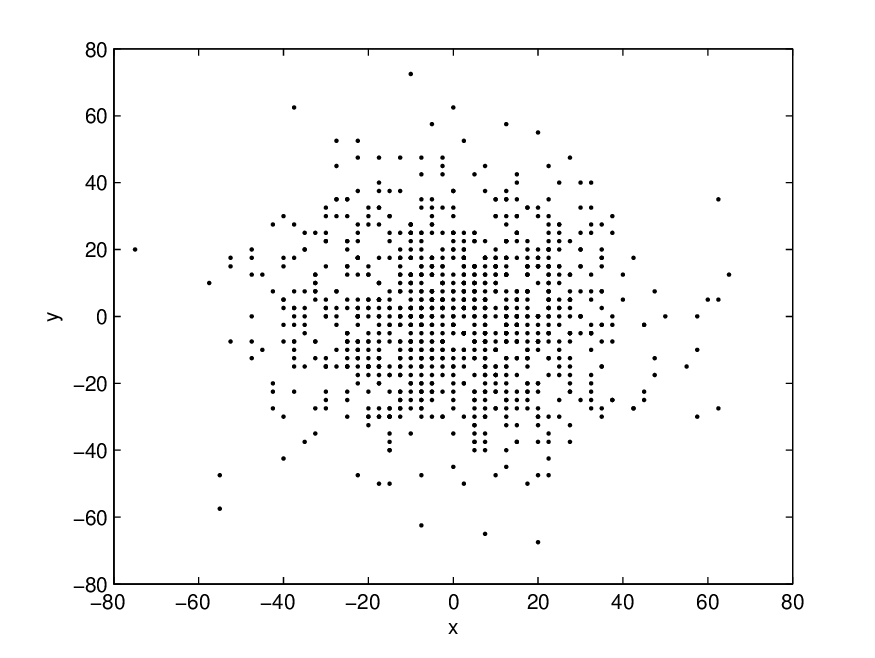}
\end{minipage}}\subfigure[]{
\begin{minipage}[t]{0.5\textwidth}\label{fig:d2}
\centering
\includegraphics[scale=0.4]{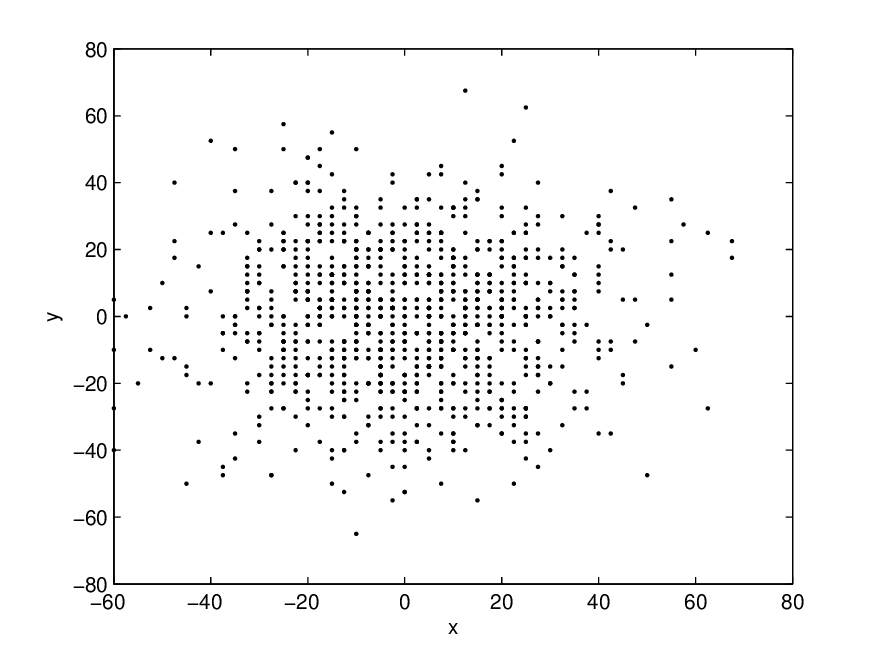}
\end{minipage}}
\caption{ Diffusion position of 1000 particles at time 5(a), 150(b), 450(c), 600(d), with $H=0.7$ and $\lambda=0.01$.}\label{fig:3}
\end{figure}

%The computer simulation indicates that $\langle x(t)\rangle=\langle x(0)\rangle, \langle y(t)\rangle=\langle y(0)\rangle$, therefore $\langle x(t)\rangle=\langle y(t)\rangle=0$, when $u(x,y,0)=\mathrm{e}^{-x^2-2y^2}$.
From \figurename \ \ref{fig:2}, one can see that $\left\langle [x(t)-\langle x(t)\rangle]^2+[y(t)-\langle y(t)\rangle]^2\right\rangle\sim t^0$ for a long time t, which implies that the stochastic process is a localization diffusion. As $0<H<1$, the larger the $H$ is, the longer the time required for MSD$\sim t^0$ and the more dispersed the particles are; the result is opposite when $\lambda$ becomes large. It is consistent with the effect of the parameter $\lambda$, which moderates the length of the jump. Fractional Brownian motion is recovered when $\lambda=0$, and its MSD is like $t^{2H}$. \figurename \ \ref{fig:3} depicts the evolution of 1000 particles when the time is 5, 150, 450, 600, respectively. \figurename \ \ref{fig:2} and \figurename \ \ref{fig:3} show that most of the particles diffuse within the bounded domain, and its size is related to $H$ and $\lambda$.

The complete algorithm is shown in Algorithm 1. First, we choose the time mesh through nonuniform time stepsizes. Next, Module 1 computes $\widetilde{U}_{M,N}^K$ by the scheme \eqref{eq:6} or \eqref{eq:9}.
%Once all $\widetilde{U}_{M,N}^K$ are known, any element of $U_{M,N}^K$ can be known.
From Module 2, one can get the MSD of stochastic process by expectation formula.
\begin{table}[H]
\begin{tabular}{l}
\hline
Algorithm 1 Probability density and MSD calculation \\
\hline
 request: $H, t_{\max}$\\
 1:\quad $t_1=(\tau k)^{1/{2H}}, t_2=(\tau k)^{1/H}$\\
 2:\quad if $H<0.5$\\
 3:\quad\quad  $t=t_1$\\
 4:\quad else \\
 5:\quad\quad $t_{\max}=(\tau k_1)^{1/{2H}}=(\tau k_2)^{1/H}$\\
 6:\quad\quad  $t=[t_1(1:k_1),t_2(k_2+1:\mathrm{end})]$\\
  \hline
  Module 1-Running  of scheme \eqref{eq:6} or \eqref{eq:9}\\
  \hline
  request: $C(t_k), C_1(t_k), C_2(t_k), U_{M,N}^0$\\
 7:\quad $\widetilde{U}_{M,N}^0=\mathrm{reshape}\left(U_{M,N}^0, (M+1)(N +1), 1\right)$\\
 8:\quad if $H<0.5$\\
 9:\quad\quad $\widetilde{U}_{M,N}^{k+1}=C^{-1}(t_k)\widetilde{U}_{M,N}^{k}$\\
 10:\quad else \\
 11:\quad while $k<k_1$\\
 12:\quad\quad $\widetilde{U}_{M,N}^{k+1}=C_1^{-1}(t_k)\widetilde{U}_{M,N}^{k}$\\
 13:\quad while $k>k_1$\\
 14:\quad\quad $\widetilde{U}_{M,N}^{k+1}=C_2^{-1}(t_k)\widetilde{U}_{M,N}^{k}$\\
15:\quad All $U_{M,N}^K$ are known \\
 \hline
 Module 2-Calculating MSD \\
 \hline
 16:\quad $\mathrm{Pr}_{M,N}^k=\frac{U_{M,N}^k}{\sum\limits_{m=0}^{M}\sum\limits_{n=0}^{N}U_{m,n}^k}$           \\
 17:\quad $\mathrm{MSD}=\sum\limits_{m=0}^{M}\sum\limits_{n=0}^{N}\left(x_{m,n}^2+y_{m,n}^2\right)\mathrm{Pr}_{m,n}^k$\\
 \hline
\end{tabular}
\end{table}
\section{Conclusion}

Anomalous diffusion is widely observed in the nature world, the types of which are abundant, and  the mechanisms of different types of anomalous diffusions sometimes are fundamentally different. This paper focuses on providing the numerical methods for the Fokker-Planck equation governing the PDF of the tfBm, and simulating the corresponding dynamics. The main challenges come from the variable coefficient of the model, and even its singularity at the starting point $t=0$. By introducing the nonuniform time stepsizes, the efficient numerical schemes are designed, and the numerical stability and convergence are theoretically proved. The simulation results, using the proposed schemes, further reveal the dynamics of the localization diffusion of tfBm.

\section*{Acknowledgements}
This work was supported by the National Natural Science Foundation of China under grant
no. 11671182.

%In this paper, the numerical schemes of the two-dimensional Fokker-Planck equation for tempered fractional Brownian motion have been discussed. When $0<H<0.5$, $2Ht^{2H-1}\partial t$ instead of $\partial \left(t^{2H}\right)$ to avoid the singularity, using nonuniform time stepsizes $t_k=(\tau k)^{1/{2H}}$ improves computational efficiency and ensures accuracy of numerical solution. For $0.5<H<1$, the time span dependent nonuniform timesteps method improves computational efficiency and the accuracy of numerical solution. These techniques which can also be applied to other numerical methods to solve the Fokker-Planck equation of tfBm do not change the original numerical method of the convergence order.
%
%The MSD of the stochastic process is obtained by numerical methods in the two-dimension space. The numerical result demonstrate that stochastic process corresponding to the problem \eqref{eq:2} is the localization diffusion process which is consistent with the characteristic of tfBm.
\appendix
\section{Numerical stability}
For $0<H<0.5$, the Fourier series of $u^k(x,y)$ is
\begin{equation}\label{eq:A.1}
%\begin{array}{ll}
u^k(x,y)=\sum\limits_{p_1=-\infty}^{+\infty}\sum\limits_{p_2=-\infty}^{+\infty}
\hat{u}_{p_1,p_2}^k\exp\left(i\frac{2p_1\pi x}{L}+i\frac{2p_2\pi y}{L'}\right),
%\end{array}
\end{equation}
where
$$
\hat{u}_{p_1,p_2}^k=\frac{1}{LL'}\int_{0}^{L}\int_{0}^{L' } u^k(x,y)\exp\left(-i\frac{2p_1\pi x}{L}-i\frac{2p_2\pi y}{L'}\right)\ud x \ud y\quad p_1,p_2=0,\pm1,{}\cdots.
$$
There exists Parseval equation
\begin{displaymath}
\left\|u^k(x,y)\right\|_{L^2}^2=LL'\sum_{p_1=-\infty}^{+\infty}
\sum_{p_2=-\infty}^{+\infty}\left|\hat{u}_{p_1,p_2}^k\right|^2.
\end{displaymath}
From \eqref{eq:6}, we get
\begin{equation}\label{eq:A.2}
\begin{array}{l}
u^{k+1}\left(x+x_m,y+y_n\right)-u^k\left(x+x_m,y+y_n\right)\\
=\frac{r}{h^2}\delta_x^2u^{k+1}\left(x+x_m,y+y_n\right)^{k+1}+
\frac{r}{l^2}\delta_y^2u^{k+1}\left(x+x_m,y+y_n\right).
\end{array}
\end{equation}
Substituting Eq. \eqref{eq:A.1} into Eq. \eqref{eq:A.2} leads to
\begin{equation}\label{eq:A.3}
\begin{array}{l}
\sum\limits_{p_1=-\infty}^{+\infty}\sum\limits_{p_2=-\infty}^{+\infty}
\hat{u}_{p_1,p_2}^kQ\left(p_1,p_2\right)\\ =\sum\limits_{p_1=-\infty}^{+\infty}\sum\limits_{p_2=-\infty}^{+\infty}
\hat{u}_{p_1,p_2}^{k+1}Q\left(p_1,p_2\right)
\bigg\{\left(1+\frac{2r}{h^2}+\frac{2r}{l^2}\right)\\
~~~~-\frac{r}{h^2}
\left[\exp\left(i\frac{2p_1\pi h}{L}\right)+\exp\left(-i\frac{2p_1\pi h}{L}\right)\right]-\frac{r}{l^2}\left[\exp\left(i\frac{2p_2\pi l}{L'}\right)+\exp\left(-i\frac{2p_2\pi l}{L'}\right)\right]\bigg\},
\end{array}
\end{equation}
where
\begin{displaymath}
Q\left(p_1,p_2\right)=\exp\left(i\frac{2p_1\pi x}{L}+i\frac{2p_2\pi y}{L'}\right)\exp\left(i\frac{2p_1\pi m h}{L}+i\frac{2p_2\pi n l}{L'}\right).
\end{displaymath}
Since the two sides of Eq. \eqref{eq:A.3} are the Fourier series, we have
\begin{equation}\label{eq:A.4}
\hat{u}_{p_1,p_2}^{k+1}=G_1\left(p_1h,p_2l\right)\hat{u}_{p_1,p_2}^k,
\end{equation}
where
\begin{displaymath}
G_1\left(p_1h,p_2l\right)=\frac{1}{1+\frac{2r}{h^2}\left(1-\cos\frac{2p_1\pi h}{L}\right)+\frac{2r}{l^2}\left(1-\cos\frac{2p_2\pi l}{L'}\right)}.
\end{displaymath}
This implies that
\begin{displaymath}
0\le G_1\left(p_1h,p_2l\right)\le1.
\end{displaymath}
Combining Parseval equation and Eq. \eqref{eq:A.4} results in
{\setlength\arraycolsep{2pt}
\begin{eqnarray*}
\left\|u^k(x,y)\right\|_{L^2}^2
&=&LL'\sum_{p_1=-\infty}^{+\infty}\sum_{p_2=-\infty}^{+\infty}
\left|\hat{u}_{p_1,p_2}^k\right|^2\\
&<&\left\|u^0(x,y)\right\|_{L^2}^2.
\end{eqnarray*}}

As $0.5<H<1$, for $t\ge t_{\max}=t_{k_1}$, using the same process, we have
\begin{equation}\label{eq:A.5}
\hat{u}_{p_1,p_2}^{k+1}=G_2\left(p_1h,p_2l\right)\hat{u}_{p_1,p_2}^k,
\end{equation}
where
\begin{displaymath}
G_2\left(p_1h,p_2l\right)=\frac{1}{1+\frac{2r_1}{h^2}\left(1-\cos\frac{2p_1\pi h}{L}\right)+\frac{2r_1}{l^2}\left(1-\cos\frac{2p_2\pi h}{L'}\right)}
\end{displaymath}
and
\begin{displaymath}
r_1=\frac{\Gamma(H+1/2)}{H\sqrt{\pi}(2\lambda)^H}
\left[\lambda t_{k+1}K_{H-1}(\lambda t_{k+1})\tau\right].
\end{displaymath}
For $k\le k_1$, with the proof being completely the same as the case that $0<H<0.5$, there exists
\begin{equation}\label{eq:A.6}
\left\|u^k(x,y)\right\|_{L^2}^2<\left\|u^0(x,y)\right\|_{L^2}^2.
\end{equation}
For $k>k_1$, combining \eqref{eq:A.5} and \eqref{eq:A.6} leads to
{\setlength\arraycolsep{2pt}
\begin{eqnarray*}
\left\|u^k(x,y)\right\|_{L^2}^2
&=&LL'\sum_{p_1=-\infty}^{+\infty}\sum_{p_2=-\infty}^{+\infty}
\left[G_2\left(p_1h,p_2l\right)\right]^{2k-2k_1}\left|\hat{u}_{p_1,p_2}^{k_1}\right|^2\\
&<&\left\|u^0(x,y)\right\|_{L^2}^2.
\end{eqnarray*}}
\section{Convergence}
We use notations
\begin{equation*}
\begin{array}{l}
Lu(x,y,t)=\frac{\partial u(x,y,t)}{\partial \left(t^{2H}\right)}-\frac{\Gamma(H+1/2)\lambda t^{1-H}K_{H-1}(\lambda t)}{2H\sqrt{\pi}(2\lambda)^H}
\left[\frac{\partial^2}{\partial x^2}+\frac{\partial^2}{\partial y^2}\right]u(x,y,t),\\
L^{(1)}u_{m,n}^k=\frac{u_{m,n}^{k+1}-u_{m,n}^k}{\bigtriangleup \left(t_k^{2H}\right)}-\frac{\Gamma(H+1/2)\lambda t_{k+1}^{1-H}K_{H-1}(\lambda t_{k+1})}{2H\sqrt{\pi}(2\lambda)^H}\left[\frac{\delta_x^2}{h^2}+\frac{\delta_y^2}{l^2}\right]u_{m,n}^{k+1}.
\end{array}
\end{equation*}  %$t^{2H}=t_k^{2H}$
As $0<H<0.5$, performing the Taylor expansion at $t_k^{2H}$, there exist
\begin{equation}\label{eq:B.1}
\begin{array}{l}
\frac{\bigtriangleup t_k}{\bigtriangleup\left(t_k^{2H}\right)}
=\frac{t_{k+1}-t_k}{\tau}\\
=\frac{t_k^{1-2H}}{2H}-\frac{(1-2H)t_k^{1-4H}}{8H^2}\tau+O\left(\tau^2\right)
\end{array}
\end{equation}
and
\begin{equation}\label{eq:B.2}
\frac{\left(\bigtriangleup t_k\right)^2}{\bigtriangleup \left(t_k^{2H}\right)}=\frac{t_k^{2-4H}}{4H^2}\tau+O\left(\tau^2\right).
\end{equation}
Letting
$
R_{m,n}^k=L^{(1)}u_{m,n}^k-[Lu(x,y,t)]_{m,n}^k,
$
and using Eq. \eqref{eq:B.1} and \eqref{eq:B.2} lead to
\begin{equation}\label{eq:B.3}
\begin{array}{l}
R_{m,n}^k=-\frac{(1-2H)t_k^{1-4H}}{8H^2}\tau\left(\frac{\partial u(x,y,t)}{\partial t}\right)_{m,n}^k-\frac{t_k^{2-4H}}{8H^2}\tau\left(\frac{\partial^2 u(x,y,t)}{\partial t^2}\right)_{m,n}^k+O\left(\tau^2+h^2+l^2\right)\\
=O\left(\tau+h^2+l^2\right)
\end{array}.
\end{equation}
For $e_{m,n}^k=u_{m,n}^k-u(x_m,y_n,t_k)$, from
Eqs.\ \eqref{eq:2}, \eqref{eq:6}, and \eqref{eq:B.3}, we have
\begin{equation*}
\begin{array}{l}
e^{k+1}\left(x+x_m,y+y_n\right)-e^k\left(x+x_m,y+y_n\right)\\
=r\left[\frac{\delta_x^2 }{h^2}+\frac{\delta_y^2 }{l^2}\right]e^{k+1}\left(x+x_m,y+y_n\right)+\tau R^k\left(x+x_m,y+y_n\right).
\end{array}
\end{equation*}
Following the proof process of numerical stability and using the expansion similar to \eqref{eq:A.3}, there exists
%Applying the proof steps of \eqref{eq:A.3}, we have
\begin{displaymath}
\left\|e^{k+1}\right\|_{L^2}^2<\left\|e^k+\tau R^k\right\|_{L^2}^2,
\end{displaymath}
leading to
{\setlength\arraycolsep{2pt}
\begin{eqnarray*}
\left\|e^k\right\|_{L^2}&<&\left\|e^{k-1}\right\|_{L^2}+\tau \left\|R^{k-1}\right\|_{L^2}\\ &\le&\left\|e^0\right\|_{L^2}+k\tau \max_{0\le i\le k}\left\|R^i\right\|_{L^2}\\ &\le&t_k^{2H}\max_{0\le i\le k}\left\|R^i\right\|_{L^2}\\&=&O\left(\tau+h^2+l^2\right).
\end{eqnarray*}}

For $0.5<H<1$, when $t>t_{\max}$, by Taylor expansion at $t_k^H$, we have
\begin{equation}\label{eq:B.4}
\frac{\bigtriangleup t_k}{\bigtriangleup \left(t_k^{H}\right)}
=\frac{t_k^{1-H}}{H}-\frac{(1-H)t_k^{1-2H}}{2H^2}\tau+O\left(\tau^2\right)
\end{equation}
and
\begin{equation}\label{eq:B.5}
\frac{\left(\bigtriangleup t_k\right)^2}{\bigtriangleup \left(t_k^{H}\right)}=\frac{t_k^{2-2H}}{H^2}\tau+O\left(\tau^2\right),
\end{equation}
%Eq. \eqref{eq:B.4}  and \eqref{eq:B.5}  show that
which implies that
\begin{equation*}
\begin{array}{l}
R_{m,n}^k=-\frac{(1-H)t_k^{1-2H}}{2H^2}\tau\left(\frac{\partial u(x,y,t)}{\partial t}\right)_{m,n}^k-\frac{t_k^{2-2H}}{2H^2}\tau\left(\frac{\partial^2 u(x,y,t)}{\partial t^2}\right)_{m,n}^k+O\left(\tau^2+h^2+l^2\right)\\
=O\left(\tau+h^2+l^2\right).
\end{array}
\end{equation*}
For $k\le k_1$,
\begin{equation}\label{eq:B.6}
\left\|e^k\right\|_{L^2}<\left\|e^0\right\|_{L^2}+t_k^{2H}\max_{0\le i\le k}\left\|R^i\right\|_{L^2}=O\left(\tau+h^2+l^2\right).
\end{equation}
When $k\ge k_1$, combining \eqref{eq:B.6} leads to
{\setlength\arraycolsep{2pt}
\begin{eqnarray*}
\left\|e^k\right\|_{L^2}&<&\left\|e^{k_1}\right\|_{L^2}+\tau (k-k_1) \left\|R^{k-1}\right\|_{L^2}\\ &\le&\left\|e^0\right\|_{L^2}+t_{k_1}^{2H} \max_{0\le i\le k_1}\left\|R^i\right\|_{L^2}+t_k^{H} \max_{k_1\le i\le k}\left\|R^i\right\|_{L^2}\\ &\le&\left(t_k^{H}+t_{\max}^{2H}\right)\max_{0\le i\le k}\left\|R^i\right\|_{L^2}\\&=&O\left(\tau+h^2+l^2\right).
\end{eqnarray*}}
\label{}



%% References


%\bibliographystyle{elsarticle-num}
%\bibliography{references}



\end{document}